\documentclass[11pt,reqno]{amsart}
\usepackage{hyperref}
\usepackage{amsmath,amsfonts,amssymb,amsthm}
\usepackage[latin1]{inputenc}
\usepackage{verbatim}
\usepackage{color}

\usepackage{amsthm}
\theoremstyle{rmk} 
\newtheorem{rmk}{Remark}[section] 

\newtheorem{theorem}{Theorem}[section]
\newtheorem{lemma}[theorem]{Lemma}
\newtheorem{prop}[theorem]{Proposition}
\newtheorem{cor}[theorem]{Corollary}

\theoremstyle{defn}

\newcommand{\sect}{\vspace{3mm} \setcounter{equation}{0} \setcounter{defn}{0} \section}

\newcommand{\w}[1]{\langle {#1} \rangle}

\newcommand{\pf}{\noindent {\bf Proof. \hspace{2mm}}}
\newcommand{\ef}{ \hfill $ \Box $ \vskip 3mm}
\newcommand{\be}{\begin{equation}}
\newcommand{\ee}{\end{equation}}
\newcommand{\bea}{\begin{eqnarray}}
\newcommand{\eea}{\end{eqnarray}}
\newcommand{\beas}{\begin{eqnarray*}}
\newcommand{\eeas}{\end{eqnarray*}}
\newcommand{\ep}{{\epsilon}}
\newcommand{\f}{\frac}
\newcommand{\na}{\nabla}

\newcommand{\bC}{{\mathbb C}}
\newcommand{\bR}{{\mathbb R}}
\newcommand{\bN}{{\mathbb N}}

\newcommand{\vF}{{\mathcal F}}

\newcommand{\vL}{{\mathcal L}}

\newcommand{\vS}{{\mathcal S}}

\newcommand{\gM}{{\mathfrak M}}

\newcommand{\gm}{{\mathfrak m}}

\newcommand{\norm}[1]{\|#1\|}\newcommand{\ip}[2]{\langle #1,#2\rangle}
\newcommand{\Cc}{C_c^\infty}

\def\p{\partial}

\def\f{\frac}

\def\na{\nabla}

\def\al{\alpha}

\def\g{\gamma}

\def\i{\infty}

\usepackage{xcolor}

\begin{document}

\title[The Kramers-Fokker-Planck equation]{Global-in-time $L^p-L^q$ estimates for solutions of the Kramers-Fokker-Planck equation}

\author{Xue Ping  WANG \and Lu ZHU}

\date{\today}

\address{Laboratoire de Math\'ematiques Jean Leray\\
UMR CNRS 6629\\
Universit\'{e} de Nantes \\
44322 Nantes Cedex 3  France \\
E-mail: xue-ping.wang@univ-nantes.fr}

\address{College of Science\\
Hohai University\\
210024 Nanjing China\\
E-mail: zhulu@hhu.edu.cn}

\subjclass[2000]{35J10, 35P15, 47A55}
\keywords{Global-in-time estimates, nonselfadjoint operators, kinetic equation, Kramers-Fokker-Planck operator}
\begin{abstract}

In this work, we prove an optimal global-in-time $L^p-L^q$ estimate  for solutions to the Kramers-Fokker-Planck equation with short range potential in dimension three. Our result shows that the decay rate as $t\to +\infty$ is the same as the heat equation in $x$-variables
and the divergence rate as $t\to 0_+$ is related to the sub-ellipticity
with loss of one third derivatives  of the Kramers-Fokker-Planck operator.
\end{abstract}

\maketitle

\sect{Introduction}

The Kramers-Fokker-Planck equation  is the evolution equation for the distribution functions describing  the Brownian motion of particles in an external field:
\be \label{kfp1} \frac{\partial W}{\partial t} = \left( -v\cdot \nabla_x  + \nabla_v \cdot(\gamma v  -\f{F(x) }{m}) + \f{\gamma k T }{m} \Delta_v \right)W,
\ee
where  $F(x) = - m \nabla V(x)$ is the external force and $W=W(t; x,v)$ is the  distribution function of particles for $x, v \in \bR^n$ and $t > 0$. In this equation,  $x$ and $v$  represent  the position and velocity variables of particles, $m$ the mass, $k$ the Boltzmann constant, $\gamma$ the friction coefficient and $T$ the temperature of the media.
This equation, called  the Kramers equation in the book of H. Risken \cite{risc}, was initially derived and used  by H. A. Kramers \cite{kr} to describe  kinetics of chemical reaction. Later on it turned out that it had more general applicability to different fields such as supersonic conductors, Josephson tunneling  junction and relaxation of dipoles.
Equation {(\ref{kfp1})}, also often called the Fokker-Planck equation, is in fact a special case of the more general  Fokker-Planck equation (\cite{risc}) or the Kolmogorov forward equation for continuous-time diffusion processes (\cite{kol}).\\
\\

After appropriate normalisation of physical constants and change of the unknown,  the KFP equation can be written into the form
\be\label{equation} \p_t u(t; x,v)+P u(t; x,v)=0,\ (x,v)\in\bR^n\times\bR^n,   t >0,
\ee with initial data
\be\label{initial} u(0; x,v)=u_0(x,v), \ee
where $P$ is the KFP  operator defined by
\be\label{operator}
 P=-\Delta_v+\f{1}{4}|v|^2-\f{n}{2} + v\cdot\na_x-\na V(x)\cdot\na_v. \ee
The main result of this work is proved for short-range potentials, that is,   $V(x)$ is  a real-valued $C^1$ function verifying
\be \label{ass1} |V(x)| + \w{x}|\nabla V(x)| \le C\w{x}^{-\rho},  \quad x\in \bR^n, \ee
for some $\rho > 1$. Here $\w{x} = (1 + |x|^2)^{1/2}$.
Remark that $V(x)$ is determined up to an additive constant. (\ref{ass1})  implies that  this constant is
chosen such that
 \[
\lim_{|x|\to \infty } V(x)= 0
\]
which can be interpreted as a normalization condition for $V(x)$.
Let $\gm$ be the function defined by
\be\label{Maxwellian} \gm(x,v)= \f{1}{(2\pi)^{\f n 4}}e^{-\f 1 2 (\f{v^2}{2} + V(x))}. \ee
Then $\gM=\gm^2$ is the Maxwellian (\cite{risc}) and $\gm$ verifies the stationary KFP equation
\be P\gm =0 \quad \mbox{ in } \bR^{2n}_{x,v}.
\ee

The large-time asymptotics of the solution to the KFP equation is mostly motivated by mathematical analysis of trend to  equilibrium in statistical physics and is studied by many authors for confining potentials. See,  for example, \cite{dv,hln,hrn,hhs}. In these works,  the potential $V(x)$ is supposed to be confining so that the spectrum of $P$ is discrete.  The typical result is  return to the equilibrium with exponential rate: $\exists \sigma >0$ such that
\be
\label{eq1.6}
u(t)= \w{\gm, u_0} \gm + O(e^{-\sigma t}), \quad t\to +\infty,
\ee
 where $V(x)$ is assumed to be normalized by
\[
\int_{\bR^n} e^{-V(x)} dx =1.
\]
In \cite{lz1}, sub-exponential convergence rate is obtained for weakly confining potential.  For quickly decreasing potentials (or more precisely,  for quickly decreasing $|\nabla V(x)|$), it is shown in \cite{nw} for $n=1$ and  \cite{w3} for $n=3$ that
\be
u(t)= \f{1}{(4\pi t)^{\f n 2}}\left(\w{\gm, u_0} \gm + O(t^{-\ep})\right), \quad t\to +\infty,
\ee
in weighted $L^2$-spaces with weight in $x$-variables. \\

In this work, we study $L^p-L^q$ estimates of the solution $u(t)$ of equation \eqref{equation}.
Here
\[
L^p = L^p(\bR^{2n}_{x,v}; dx dv)
\] is equipped with the natural norm. For $f \in L^p$ and $T$ bounded linear operator from $L^p$ to $L^q$, we denote :
\be
\|f\|_p = \|f\|_{L^p}, \quad  \|T\|_{p \to q } =\|T\|_{\vL(L^p, L^q)}.
 \ee
By an abuse of notation, for a closed linear operator $T$ in $L^2$ with $C_0^\infty(\bR^{2n})$ as a core and for $p\in[1, \i[$,  we still denote by the same letter  $T$ its minimal closed extension in $L^p$ ({\it i.e.}, the closure in $L^p$ of the restriction of $T$ to $C_0^\infty(\bR^{2n})$).
Similarly, the notation  $e^{-tP} : L^p \to L^q$ means that the restriction of $e^{-tP}$ on $C_0^\infty$ extends to a map from $L^p$ to $L^q$.
Under fairly general condition, $e^{-tP}$ is a strongly continuous positivity preserving  semigroup in $L^p$.
We shall prove that if $V\in C^1$ with $\|V\|_\infty, \|\nabla_xV\|_\infty<\infty$,  then $e^{-tP}$ is uniformly bounded on $L^p(\bR^{2n})$ for  $p \in[1, \i[$, the infinitesimal generator of this semigroup being $-P$, where $P$ is the minimal closed extension of the KFP operator in $L^p$:
\[
e^{-tP} = \overline{ \left( e^{-tP}|_{C_0^\i} \right) }|_{L^p}  = e^{-t \,\overline{( P|_{C_0^\i} )}|_{L^p}}
\]
Since for $1 \le p <\i$,
\[
\overline{ \left( e^{-tP}|_{C_0^\i} \right) }|_{L^p}  = e^{-t \,\overline{( P|_{C_0^\i} )}|_{L^p}},
\]
 our notation is consistent. When $p=\infty$,  $e^{-tP}$ on $L^\infty$ is defined as the dual of $e^{-tP^*}$ on $L^1$. The main result of this article is the following

\begin{theorem}\label{th1.1} Let $n=3$ and condition (\ref{ass1}) be satisfied with $\rho >1$.  For $1 \le p < q \le \i$, there exists some constant $C>0$ such that
\be\label{e1.12} \|e^{-tP}\|_{p \to q} \le \f{C}{(\gamma(t))^{\f 3 {2 p}(1-\f p q)}}, \quad t\in ]0, \i[,
\ee
 where $\gamma(t) = \sigma(t) \theta(t) $ with
\be \label{e1.13}
  \quad \sigma(t) = t - 2 \coth (t) + 2 \text{\rm cosech}(t),  \quad  \theta(t) = 4\pi  e^{- t } \sinh (t).
\ee
\end{theorem}

The function $\gamma(t)$ appears in the explicit formula  of  fundamental solution for the free KFP equation (see Section 2) and behaves like: $\gamma (t) \sim 2\pi t $ as $t\to\i$ and
{$\gamma(t)\sim \f\pi3 t^4$ as $t \to 0_+$.}
The two factors of $\gamma(t)$ have different meanings. $\theta(t)$  arises from  the semigroup generated by the harmonic oscillator
\[
H = \Re P = -\Delta_v+\f{1}{4}|v|^2-\f{n}{2}
\]
in $L^p(\bR^3_v)$. For $\sigma(t)$,  remark   that
\be
\sigma(t) \sim t,  \mbox{ as } t \to +\infty; \quad \sigma(t) \sim \f{t^3}{12},  \mbox{ as } t \to 0_+.
\ee
For $p=1, q =\i$,
\[
 (\sigma(t))^{-\f 3 {2 }} \sim \f{C_1}{t^{\f {3} {2}}}, \quad t\to \i; \quad  (\sigma(t))^{-\f 3  2} \sim \f{C}{t^{\f 9  2}}, \quad t \to 0_+.
\]
One sees that the term $(\sigma(t))^{-\f 3 {2 }} $ is of the same order  as that of the heat semigroup $e^{t\Delta_x}$ as map from $L^1(\bR^3)$ to $L^\i(\bR^3)$ as $t\to \i$ and of  the same order as that of $e^{-t|D_x|^{\f 2 3}}$ as $t\to 0_+$. This may be explained by the fact that at low energies, the KFP operator $P$  behaves like a Witten Laplacian (\cite{hln,lwx}), while globally it is  sub-elliptic in $x$ with the loss of $\f 1 3$ derivatives.  \\

To prove Theorem \ref{th1.1}, we  first study   the semigroup  $e^{-tP_0}$ in $L^p-L^q$ setting, where
\[
P_0 = -\Delta_v+\f{1}{4}|v|^2-\f{n}{2} + v\cdot\na_x.
\]
Then we consider $P$ as perturbation of $P_0$ and use Duhamel's formula to prove (\ref{e1.12}).
The short-time estimate for $e^{-tP}$ can be easily obtained (see Theorem \ref{short}) and is valid for $n \ge 1$ and for any real-valued $V\in C^1$ with $\na V$ bounded.
The proof of (\ref{e1.12}) for $t $ large is based  on a result of  time-decay of $e^{-tP}$ in weighted $L^2$ spaces obtained in \cite{w3}.
\\

In this work, we often use an argument of duality which is based on  the relation
\be \label{Jrel}
P^* = JPJ
\ee
 in $L^2$, where $J$ is the reflection in $v$ variable: $Jf (x, v)  = f(x,-v)$.  If one has some estimates for $P$ or {$e^{-tP}$}  in $L^p$, one can often use the  duality between $L^p$ and $L^q$,  $p^{-1} + q^{-1} =1$,  to affirm  that the same statements are true for $P^*$ or $(e^{-tP})^*$ in $L^q$.  Since $J$ preserves any $L^p$ norm, the same estimates hold true for $P$ or $e^{-tP}$ in $L^q$.
  \\

The remaining part of this work is organized as follows.  In Section 2, we establish an explicit useful formula for the fundamental solution of the free KFP operator $P_0$. Global-in-time $L^p-L^q$ estimates are obtained for $e^{-tP_0}$ in Section 3.  Theorem \ref{th1.1} is proved in Section 4.

\sect{Fundamental solution of the free KFP equation}

In this Section, we use the method of complex deformation to calculate the fundamental solution of the free KFP equation.
Let $P_0$ be the free KFP operator:
\be\label{P_0} P_0=v\cdot\na_x-\Delta_v+\f{1}{4}|v|^2-\f{n}{2}, (x,v)\in \bR^{2n}. \ee
In $L^2$, using the partial
Fourier transform in $x$-variables, we have for $u \in D(P_0)$
\bea
P_0 u (x, v)& = &\vF_{x\rightarrow\xi}^{-1}\hat{P}_0(\xi) \hat{u}(\xi, v), \quad \mbox{ where }\\
\hat{P}_0(\xi) &=&-\Delta_v+\f{1}{4}\sum^n_{j=1}(v_j+2i\xi_j)^2-\f{n}{2}+|\xi|^2 \\
\hat{u}(\xi, v) & = & (\vF_{x\rightarrow\xi}u)(\xi, v) \triangleq \int_{\bR^n} e^{-i x\cdot\xi}u(x, v) \; dx.
\eea
Denote
\be
D(\hat{P}_0) =  \{f \in L^2(\bR^{2n}_{\xi, v}); \hat{P}_0(\xi)f   \in L^2(\bR^{2n}_{\xi, v})\}.
\ee
Then $\hat{P}_0  \triangleq\vF_{x\rightarrow\xi} P_0 \vF_{x\rightarrow\xi}^{-1}$ is a direct integral of the  family of complex harmonic operators $\{\hat{P}_0(\xi); \xi \in \bR^n \}$.
$\{\hat{P}_0(\xi), \xi\in \bR^n\}$ is  a holomorphic family of type $(A)$ in sense of Kato with constant domain
$D= D(-\Delta_v+\f{v^2}{4})$ in $L^2(\bR^n_v)$.
Let $F_j(s)=(-1)^je^{\f{s^2}{2}}\f{d^j}{ds^j}e^{-\f{s^2}{2}}, j \in \bN,$ be the Hermite polynomials and
$$
\varphi_j(s)=(j!\sqrt{2\pi})^{-\f{1}{2}}e^{-\f{s^2}{4}}F_j(s)
$$
the normalized Hermite functions.  For $\xi \in \bR^n$ and $\alpha=(\alpha_1, \alpha_2, \cdots, \alpha_n) \in \bN^n$, define
\be \psi_\alpha(v) = \prod_{j=1}^n\varphi_{\alpha_j}(v_j) \mbox{ and } \psi_\alpha^\xi(v) = \psi_\alpha(v + 2i \xi). \ee
Then
\be \label{e2.7}
\hat{P}_0(\xi) \psi_\al^\xi = (|\al| + |\xi|^2) \psi_\al ^\xi.
\ee
For $\alpha, \beta \in \bN^n$,  $\xi \to \w{\psi_\alpha^\xi, \psi_\beta^{-\xi}} $ extends to an entire function for $\xi \in \bC$ and is constant on $i\bR$. Therefore $\w{\psi_\alpha^\xi, \psi_\beta^{-\xi}}$ is  constant for $\xi \in \bC$ and one has
\be \label{basis} \w{\psi_\alpha^\xi, \psi_\beta^{-\xi}} = \delta_{\alpha\beta} = \left\{\begin{array}{ll}
$1$, & \hbox{$\alpha =\beta$,} \\
$0$, & \hbox{$\alpha \neq \beta$.}
\end{array}
\right., \quad \forall  \alpha, \beta \in \bN^n, \xi \in \bR^n. \ee
$e^{-tP_0}$ is a {contraction} semigroup in $L^2(\bR^{2n}_{x,v})$ {($P_0$ is accretive: $\Re\w{P_0u,u}=\w{Hu,u}\ge0$ for $u \in C_0^\infty$)}. Its distributional kernel  can be explicitly computed, using  Mehler's formula for harmonic oscillator (\cite{meh})(see also \cite{boch}), where this fundamental solution is calculated with different method and expressed in slightly different way.  Recall (\cite{cfks}) that the heat kernel of the $n$-dimensional harmonic oscillator $-\Delta + x^2$ is given by
\be \label{mehler}
E(x,y; t) = \f{1}{(2 \pi \sinh (2t))^{\f n 2}}\exp\left({  -\f{\coth (2t)}{2} (x^2 + y^2) + {\text{cosech} (2t)}
x\cdot y}\right), \quad t>0.
\ee

\begin{lemma} \label{lem2.1} Let $n \ge 1$. The distributional kernel of $e^{-tP_0}$ is given by
\be \label{Fvt}
F(x,v, x',v';t) =  \f{1}{(4\pi \sigma(t))^{\f n2}} \exp\left({- \f{1}{4\sigma(t)} |x-x'-\omega(t) (v+v')|^2}\right) K(v,v';t). \ee
where
\bea
K(v,v'; t) &=& \f{1}{(4 \pi \sinh (t))^{\f n 2}}\exp\left({\f {nt} 2  -\f{\coth (t)}{4} (|v|^2 + |v'|^2) +
\f{\text{\rm cosech} (t)}{2} v \cdot v'}\right)  \label{Kvt} \\[2mm]
\omega(t) &= &  \coth(t) -{\text{\rm cosech} (t)} \nonumber
\\[2mm]
 \sigma(t) &=& t - 2 \coth (t) + 2 \text{\rm cosech}(t). \nonumber
\eea
\end{lemma}
\pf Since the $n$-dimensional free KFP operator $P_0$  is a direct
 sum of $n$ one dimensional operators, it suffices to prove the lemma for $n=1$. Applying (\ref{mehler}) and making use of change of scale, we deduce that the  Mehler's formula for the heat kernel of the one-dimensional harmonic oscillator  $H= -\f{d^2}{dv^2} + \f 1 4 v^2 -\f 1 2$ is given by:
\be
e^{-tH} u = \int_\bR K(v,v'; t) u(v') dv', \quad t>0, u \in C_0^\infty,
\ee
 where \be
 K(v,v'; t) = \f{1}{\sqrt{4 \pi \sinh (t)}}\exp\left({\f t 2  -\f{\coth (t)}{4} (v^2 + v'^2) + \f{\text{cosech} (t)}{2}
v v'}\right), \ee
which is an entire function in $v$ and $v'$ in $\bC$.
Set
\be \tilde K(v, v', \xi; t) = e^{- |\xi|^2 t} K(v+ 2i \xi, v'+ 2i \xi;t). \ee
Since $\psi_l^\xi$ is an eigenfunction of $\hat P _0(\xi)$ associated with the eigenvalue $l + |\xi|^2$, one has
\[
e^{-t \hat P _0(\xi)} \psi_l^\xi = e^{- t(l+ |\xi|^2)}\psi_l^\xi
\]
On the other hand, one has
\[
\int_\bR  K(v + 2 i\xi, v';  t) \psi_l(v') dv' = e^{-tl}\psi_l^\xi,
\]
since the both sides are entire functions in $v \in \bC$. Using deformation of contour and the decay properties of $K(v,v';t)$, one obtains
\be
\int_\bR K(v + 2 i\xi, v'+ 2i\xi; t) \psi_l^\xi (v') dv' =  \int_\bR K(v + 2 i\xi, v';  t) \psi_l(v') dv'
\ee
for $\xi \in \bR$.
It follows that
\[
e^{-t \hat P_0(\xi)} \psi_l^\xi = \int_\bR \tilde K(v, v', \xi; t) \psi_l^\xi(v') dv' =e^{- t(l+ |\xi|^2)}\psi_l^\xi.
\]
Since the span of $\{\psi_l^\xi, l \in \bN\}$ is dense {in} $L^2(\bR_v)$, one concludes that the heat kernel of
 $\hat{P}_0(\xi)$ is equal to $\tilde K(v, v', \xi; t)$ for $t>0$.
 $\tilde K(v, v', \xi; t)$ can be written as
 \be \tilde K(v, v', \xi; t) =   K(v,v';t) \hat g(v, v', \xi; t) \ee
  where
  \be \hat g(v, v', \xi;t) = \exp\left({ -i \omega(t) (v+v')\xi -|\xi|^2 \sigma(t)}\right) \ee
  with
  \be \omega(t) =   \coth(t) -{\text{cosech} (t)}, \quad \sigma(t) = t - 2 \coth (t) + 2 \text{cosech }(t). \ee
Since $\sigma(t) >0$ for $t>0$,   the inverse Fourier transform of $\hat g$ in $\xi$ can be explicitly calculated:
\beas
g(v,v', x; t) &=&\f{1}{2\pi}\int_\bR e^{i x \xi}\hat{ g}(v,v', \xi;t) \; d\xi \\
&=&  \f{1}{2\pi}\int_\bR e^{i (x -\omega(t) (v+v')) \xi} e^{ - \sigma(t)|\xi|^2 }\; d\xi \\
& = &\f{1}{\sqrt{4\pi \sigma(t)}} \exp\left({- \f{1}{4 \sigma(t)} (x-\omega(t) (v+v'))^2}\right). \eeas
Therefore, the integral kernel of $e^{-tP_0}$ is given by
\be F(x,v,x',v';t) =  \f{1}{\sqrt{4\pi \sigma(t)}} \exp\left({- \f{1}{4 \sigma(t)} (x-x'-\omega(t) (v+v'))^2}\right) K(v,v';t). \ee \ef

The fundamental solution $F(x,v,x',v';t)$ for the free KFP equation has several nice properties. For example,  one has for $f \in C_0^\i(\bR^{2n})$,
\be \label{mean}
\int (e^{-tP_0}f) (x,v) dx  = (e^{-tH} g) (v), \quad v \in \bR^n,
\ee
where $g(v) = \int f(x, v) dx$ and  $H$ is the harmonic oscillator: $H= -\Delta_v + \f 1 4 v^2 -\f n 2$.
Note that $F>0$, and $F$ is symmetric,
\[
F(x,v,x',v';t)=F(x',v',x,v;t).
\]

\sect{Global-in-time estimates for the free KFP operator}

In this section, we give some global-in-time $L^p-L^q$ estimates for $e^{-tP_0}$ needed in the proof of Theorem \ref{th1.1}.

\begin{prop} \label{prop3.1} Let $n \ge 1$. For $t>0$, $e^{-tP_0}$ defined on $C_0^\infty(\bR^{2n})$ extends to an operator bounded from $L^1$ to $L^\i$  and  the following  estimate is true for the free KFP operator:
\be \label{e3.1}
\|e^{-tP_0}\|_{1 \to \i} \le \f{1}{(4\pi\gamma(t))^\f{n}{2}}
\ee
 for $t>0$. Here where
\[
\gamma(t) = \sigma(t) \theta(t), \quad  \theta (t)  = 4 \pi  e^{-  t  }  \sinh (t).
\]
\end{prop}
\pf   Let $f \in C_0^\infty(\bR^{2n})$.  By Lemma \ref{lem2.1}, one has
\be \label{e3.2b}
|e^{-tP_0}f (x,v)|  \le  \f{1}{(4\pi \sigma(t))^{\f n2}} \int_{\bR^{2n}} K(v, v';t) |f(x',v')| dx'dv'
\ee
which gives
\bea
|e^{-tP_0}f (\cdot,v)|_{L^\i_x}  &\le & \f{1}{(4\pi \sigma(t))^{\f n2}} \int_{\bR^{n}} K(v, v';t) \|f(\cdot,v')\|_{L^1_x}dv' \nonumber \\
 &= &  \f{1}{(4\pi \sigma(t))^{\f n2}} (e^{-tH}g) (v) \label{etH}
\eea
where $g(v') = \|f(\cdot,v')\|_{L^1_x} $, since $K(v, v',t)$ is the distributional kernel of $e^{-tH}$. From (\ref{Kvt}), it follows that
\be
\|e^{-tP_0}f\|_{\i}  \le  \f{1}{(4\pi \gamma (t))^{\f n2}} \|f\|_{1}, \quad f \in C_0^\infty(\bR^{2n}).
\ee
{(\ref{e3.1})} is derived by an argument of density.
\ef

Part of  following results may be known. We include a proof for reason of completeness.

\begin{cor} \label{cor3.2}
(a). One has
\be \label{Lp}
\| e^{-tP_0}\|_{p \to p} \le \Big(\f{2}{1+e^{-2t}}\Big)^{n|\f1p-\f12|}\le 2^{\,n|\f1p-\f12|}
\ee
for $1\le p\le \i$  and
\be \label{LpLq}
\|e^{-tP_0}\|_{p \to q} \le  \f{2^{{\f n 2|\f 1p+\f 1q-1|} }}{(4\pi \gamma(t))^{\f n {2 p}(1-\f p q)}}, \quad t>0,
\ee
for $1 \le p\le q \le \i$.

 (b).  $e^{-tP_0}$, $t\ge 0$,  is a strongly continuous positivity preserving  semigroup in $L^p$ for $1 \le p <\infty$, and $C_0^\i(\bR^{2n})$ is a core of $P_0$ in $L^p$.
\end{cor}
\pf (a). $P_0$ is closed and accretive in $L^2$. Therefore $e^{-tP_0}$, $t\ge 0$, is a strongly continuous  semigroup in $L^2$.
In particular, (\ref{Lp}) is true for $p=2$. We denote by the same symbol $e^{-tP_0}$ the operator induced on $L^{p}(\bR^{2n})$.
 By (\ref{mean}),  one has
\be
\|e^{-tP_0}f\|_{1} \le  \|e^{-tH}{g_{|f|}}\|_{{L^1_v}} \le \left( \f{2}{1+e^{-2t}}\right)^{\f n2} \|f\|_1  \le  2^{\f n 2}\|f\|_1
\ee
for $f \in L^1$, where $g_{|f|}(v)=\|f(\cdot,v)\|_{L^1_x}$. The first estimate implies  (\ref{Lp}) for $p =1$.
Since $\|e^{-tP_0}\|_{2\to 2} \le 1$, the case $p=\i$ follows by duality, and the general case by interpolation. We obtain (\ref{Lp}) for $p\in[1, \i]$.
(\ref{LpLq}) follows  from (\ref{e3.1}) and (\ref{Lp}) by interpolation. \\

(b).  Since $C_0^\infty (\bR^{2n})$ is a common  core of $P_0$ in $L^p$, $1\le p <\i$, the semigroup property of $e^{-tP_0}$ in $L^p$ follows from that of $e^{-tP_0}$
 in $L^2$.
The strongly continuity of $e^{-tP_0}$ in $L^1$ follows from Lemma \ref{lem2.1} by taking the limit $t\to 0_+$.
 The general case $1 < p <\i$ is deduced from the cases $p=1$ and $p=2$.  $e^{-tP_0}$ is positivity preserving, because its distributional kernel $F(x,v, x',v';t) $ is positive.
\ef

\begin{rmk}\label{rmk:cp} For $1\le p\le\i$ set
\be \label{cp}
c_p= n|\f 1 p -\f 1 2|.
\ee
(\ref{Lp}) shows that
\be \label{quasi0}
\|e^{-tP_0}\|_{p\to p}\le e^{c_pt}, \qquad t\ge0 ,
\ee
that is, $e^{-t(P_0+c_p)}$ is a contraction semigroup on $L^p$.
\end{rmk}

To study the full KFP operator $P$, we want to treat the $W = -\nabla V(x)\cdot\nabla_v$ as perturbation and need  some more  estimates for $e^{-tP_0}$. \\

\begin{prop} \label{prop3.3} Let $k \in \bN$ and $\alpha \in \bN^n$ with $|\alpha| =k$.  The following estimates
are true for the free KFP equation:
\be
\label{v1}\|\w{v}^k e^{-tP_0}\|_{1 \to \i}+\| \partial_v^\alpha e^{-tP_0}\|_{1 \to \i} \le \f{C}{(\gamma(t))^{\f n 2}}\left( 1 + t^{-\f{k}{2}}\right)
\ee
and  for any $p\in[1, \i]$,
\be
\label{v2}\|\w{v}^k e^{-tP_0}\|_{p \to p}+\|\partial_v^\alpha e^{-tP_0}\|_{p \to p}\le C\left(1 + t^{-\f{k}{2}}\right)
\ee
for $t>0$.
\end{prop}
\pf Remark that the distributional kernel $K(v,v', t)$ of $e^{-tH}$ satisfies the estimate
\[
0 \le K(v,v',t) \le  \f{1}{( \theta(t))^{\f n 2}} e^{- \f{\cosh^2(t) -1 }{2 \sinh(2t)} |v|^2}
\]
uniformly in $v'$ and that $\f{\cosh^2(t) -1 }{2 \sinh(2t)} \sim c t$ as $t \to 0$, $c>0$.
As in the proof of Proposition \ref{prop3.1}, one has for $f\in C_0^\infty (\bR^{2n})$
\beas
\|\w{v}^ke^{-tP_0}f(\cdot, v) \|_{L^\i_x}  &\le &   \f{1}{(4\pi \sigma(t))^{\f n2}}  \int_{\bR^{n}} \w{v}^k K(v, v',t)  \| f(\cdot, v')\|_{L^1_x} dv' \nonumber \\
& \le  &  \f{C}{(\gamma(t))^{\f n2}} (1 + t^{-\f k 2}) \|f\|_1, \quad t >0.
\eeas
This shows
\[
\|\w{v}^ke^{-tP_0}\|_{1 \to \i } \le  \f{C}{(\gamma(t))^{\f n2}} (1 + t^{-\f k 2}), \quad t >0.
\]
Similarly, one can estimate $\| \p_v^\alpha e^{-tP_0}\|_{1 \to \i }$ by evaluating $\sup_{v, v'} |\p_v^\alpha K(v, v',t)|$ for $t>0$ and $\alpha \in \bN^n$. (\ref{v1}) is proved.

In the same way, one has
\beas
\|\w{v}^ke^{-tP_0}f \|_{1}  &\le &   C  \int_{\bR^{n}} {\left( \int_{\bR^n} \w{v}^k K(v, v';t)\, dv \right)} \| f(\cdot, v')\|_{L^1_x} dv' \nonumber \\
 &\le  &  C {\left( \sup_{v'} \int \w{v}^k  K(v, v';t)\, dv\right)}  \|f\|_1 \label{etHk} \\
& \le  &  C_1 (1 + t^{-\f k 2}) \|f\|_1, \quad t >0.
\eeas

By arguments of duality and interpolation, we obtain for $p \in [1, \i]$
\[
\|\w{v}^k e^{-tP_0}\|_{p \to p} \le C\left(1 + t^{-\f{k}{2}}\right),  \quad t >0.
\]
Again using the formula of $K(v, v',t)$, one can show
\[
\|\p_v^\alpha e^{-tP_0}\|_{p \to p} \le C_\alpha\left(1 + t^{-\f{|\alpha|}{2}}\right),  \quad t >0.
\]
for any $\alpha \in \bN^n$. This proves (\ref{v2}).
\ef

As consequence of Proposition \ref{prop3.3}, one obtains  the following

\begin{cor}\label{cor3.4}
 For $1 \le p\le q \le \infty$ and for any $k \in \bN$ and $\alpha \in \bN^n$ with $|\alpha|=k$,  one has
\be
\label{v3}\|\w{v}^k e^{-tP_0}\|_{p \to q} + \|\partial_v^\alpha e^{-tP_0}\|_{p \to q}\le \f{C}{(\g(t))^{\f{n}{2p}(1-\f{p}{q})}}\left(1 + t^{-\f{k}{2}}\right),
\ee
and
\be
\label{v4}\| e^{-tP_0} \w{v}^k\|_{p \to q} + \|e^{-tP_0}\partial_v^\alpha \|_{p \to q}\le \f{C}{(\g(t))^{\f{n}{2p}(1-\f{p}{q})}}\left(1 + t^{-\f{k}{2}}\right),
\ee
for $t>0$.
\end{cor}

\sect{Global-in-time estimates for $e^{-tP}$}

Set $P= P_0+ W$ with $W = -\nabla V(x)\cdot\nabla_v$.
 Recall that $P_0$ verifies the following  global sub-elliptic estimate (\cite{hln}):
\be \label{e4.0}
 \|\w{v}^2 u\|_2 +  \|\w{D_v}^2 u\|_2  +  \|\w{D_x}^{\f 2 3} u\|_2  \le C (\|P_0u\|_2 + \|u\|_2), \quad u \in C_0^\infty(\bR^{2n}).
\ee
 Under the condition (\ref{ass1}) with $\rho > -1 $, $W$ is relatively compact perturbation of $P_0$ and $P$ is closed with $D(P)=D(P_0)$. Since
\[
e^{-tW}f(x, v) = f(x, v+ t\nabla V(x)),
\]
 $e^{-t P}$ is a strongly continuous positivity preserving semigroup on $L^2(\bR^{2n})$ satisfying a global bound
 $\|e^{-tP}\|_{p\to p} \le e^{c_pt}$, where $c_p$ is given by (\ref{cp}). This estimate is rough and can be improved in many cases (see, for example, Proposition \ref{prop:global} below). We are interested in $e^{-tP}$ when it is regarded as map from $L^p$ to $L^q$, $p<q$.

\subsection{Short-time estimates for $e^{-tP}$}

Assume  $V\in C^1$ with $\|\nabla_x V\|_\infty <\infty$.
For $f \in \vS(\bR^{2n})$ and $u(t) =e^{-tP}f$, one has
\[
u(t) = e^{-tP_0}f {+}\int_0^t e^{-(t-s)P_0} \big(\nabla V\cdot\nabla_v\big) u(s)\,ds
\]
By Corollary \ref{cor3.2} and Corollary \ref{cor3.4}, one obtains
\[
\|u(t)\|_p \le C\big(\|f\|_p + \int_0^t |t-s|^{-\f 1 2} \|u(s)\|_p ds\big)
\]
The standard fractional Gronwall argument gives that for any $T>0$, there is $C = C(T, \|\nabla V\|_\infty, p)>0$ such that
\[
\|u(t)\|_p \le C \|f\|_p, \qquad 0 \le t \le T.
\]
Therefore $e^{-tP}$ is locally bounded in $t$:
\be \label{shortbound}
\|e^{-tP}\|_{p \to p} \le C, \qquad 0 \le t \le T.
\ee

\begin{theorem}\label{short} Let  $n\ge 1$. Assume  $V\in C^1$ with $\|\nabla_x V\|_\infty <\infty$.  Then, for $1 \le p < q \le \i$, there exists $C>0$ such that
\be \label{e4.1}
\|e^{-tP}\|_{p \to q} \le \f{C}{\g(t)^{\f{n}{2p}(1 -\f pq)}}, \quad t \in ]0, 1].
\ee
\end{theorem}
\pf The proof is  based on Duhamel's formula
\be \label{e4.2}
e^{-tP} = e^{-tP_0} - \int_0^t e^{-(t-s)P_0} W e^{-sP} ds.
\ee
Set
\[
\alpha(p, q) =\f{n}{2}(\f 1 p -\f 1 q).
\]
Remark that $\gamma(t) \sim c t^4$ as $t\to 0 _+$. For $1 \le p \le p' \le 2$ such that  $\f 1 p - \f 1 {p'} < \f 1 {4n}$, one has: $4\alpha(p,p') < \f 1 2$. Since $e^{-tP}$ is uniformly bounded on $L^p$ for $t\in[0,1]$,  by (\ref{v4}), one has
\beas
\|e^{-tP}\|_{p\to p'} & \le & \|e^{-tP_0}\|_{p\to p'} +  C \int_0^t \| e^{-(t-s)P_0}\nabla_v\|_{p\to p'} \|  e^{-sP}\|_{p\to p} ds
\\
& \le &  C \left(\gamma(t)^{-\alpha(p,p')} + \int_0^t |t-s|^{-\f 1 2 - 4 \alpha(p,p') } ds \right) \\
&\le & C_1 \gamma(t)^{-\alpha(p,p')}, \quad \mbox{ for }  t\in ]0, 1].
\eeas

For $n \ge 1$, take $k=k(n)$ numbers $p_1, \cdots, p_{k+1}$ such that
\[
1=p_1 <p_2 <\cdots < {p_{k}} <p_{k+1}=2 \quad \mbox {and} \quad  \f 1 {p_j} - \f 1 {p_{j+1}} < \f 1 {4 n}.
\]
Writing $ e^{-tP}$ as $(e^{-\f t k P})^k$, one obtains
\beas
\|e^{-tP}\|_{1 \to 2} &\le &\| e^{-\f t k P}\|_{p_1\to p_2} \cdots  \| e^{-\f t k P}\|_{p_{k} \to 2} \\
 &\le & C {\gamma(\tfrac tk)}^{-\alpha(1,p_2) - \dots -\alpha (p_{k}, 2) }\\
&\le &C_k \gamma(t)^{-\alpha(1,2) }
\eeas
for  $t\in ]0, 1]$. This proves (\ref{e4.1}) for $p=1$ and $q=2$. The general case follows by duality and interpolation.
\ef

\subsection{Large-time  estimate for $e^{-tP}$}

We first prove a time-global uniform bound for $e^{-tP}$ on $L^p$. \\

\begin{prop}\label{prop:global}
Let $V\in C^1$ with $\|V\|_\infty, \|\nabla V\|_\infty<\infty$. Then for $1\le p\le\infty$, there exists {$B=B(n,p,\|V\|_\i,\|\na V\|_\i)>0$} such that
\[
\norm{e^{-tP}}_{p\to p}\le B, \qquad t \ge 0.
\]
\end{prop}
\pf $e^{-tP}$ is positivity preserving. For $0\le u_0\in\Cc$,  $u(t)=e^{-tP}u_0\ge 0$. Then $u\in L^1$ and $Pu\in L^1$. Put $f={\mathcal M}^{1/2}u$ with ${\mathcal M}=e^{-|v|^2/2-V}$.
{Since $J\mathcal M^{1/2}=\mathcal M^{1/2}$,  $P\mathcal M^{1/2}=0$ and ${\mathcal M}^{1/2}P = P^*{\mathcal M}^{1/2}$,
we have
\begin{align}
\frac{d}{dt}\int f(t)
&=\frac{d}{dt}\int {\mathcal M}^{1/2}u
=\int {\mathcal M}^{1/2}(\partial_tu)
=-\int {\mathcal M}^{1/2}(Pu)
=-\int (P^*{\mathcal M}^{1/2}) u\nonumber\\
&=-\int (JPJ{\mathcal M}^{1/2}) u
=-\int J(P{\mathcal M}^{1/2})u\nonumber\\
&=0,\nonumber
\end{align}
where $\int$ denotes integration with respect to
$(x,v)$ over $\bR^{2n}$.
It follows that $\int f(t)=\int f(0)$. Hence
\begin{equation}\label{eq:mass}
\int u(t)\,e^{-|v|^2/4}\le e^{\norm V_\infty/2}\int f(t)\le e^{\norm V_\infty/2}\int f(0)\le e^{\norm V_\infty}\norm{u_0}_1 .
\end{equation}
For $\delta,\epsilon\in(0,\frac 12]$ let $\psi(x)=e^{-\delta\langle x\rangle}$, $w(v)=e^{-\epsilon|v|^2/4}$, and $\Phi(t)=\int u(t)\psi w\,dxdv$; this is finite since $\psi w\in L^2$. The function $\psi w$ lies in $D(P^*)$ and $\Phi\in C^1$ with $\Phi'=-\ip{u}{P^*(\psi w)}$, where $P^*=H-T$, $T=v\cdot\nabla_x-\nabla V\cdot\nabla_v$, and where we have set $G=\norm{\nabla V}_\infty$. A direct computation gives
\[
H w=\Big[(1-\epsilon^2)\frac{|v|^2}{4}-(1-\epsilon)\frac n2\Big]w,\qquad T(\psi w)=w\,v\cdot\nabla\psi+\frac\epsilon2\,\psi\,(\nabla V\cdot v)\,w .
\]
Using $|\nabla\psi|\le\delta\psi$, $\delta|v|\le\frac{|v|^2}{16}+4\delta^2\le\frac{|v|^2}{16}+1$, $\frac\epsilon2G|v|\le\frac{|v|^2}{16}+G^2$, $1-\epsilon^2\ge\frac34$ and $u\ge0$,
\[
\Phi'\le-\int u\,\psi w\Big[\frac{|v|^2}{16}-\frac n2-1-G^2\Big]\le-\Phi+\Big(2+\frac n2+G^2\Big)\int_{|v|\le R_0}u ,
\]
because, with $R_0=4\big(2+\frac n2+G^2\big)^{1/2}$, the bracket is $\ge1$ for $|v|\ge R_0$ and $\ge-(1+\frac n2+G^2)$ everywhere, and $\psi w\le1$. By \eqref{eq:mass}, $\int_{|v|\le R_0}u\le e^{R_0^2/4}e^{\norm V_\infty}\norm{u_0}_1$. Since $\psi w\le1$ gives $\Phi(0)\le\norm{u_0}_1$, Gronwall gives, with $B=1+\big(2+\frac n2+G^2\big)e^{R_0^2/4}e^{\norm V_\infty}$, that $\Phi(t)\le e^{-t}\Phi(0)+(B-1)\norm{u_0}_1\le B\norm{u_0}_1$, uniformly in $\delta,\epsilon$. Letting $\delta,\epsilon\to0$ (monotone convergence) gives $\norm{u(t)}_1\le B\norm{u_0}_1$. For general $u_0$ use $|e^{-tP}u_0|\le e^{-tP}|u_0|$ and density. The case $p=\infty$ follows by duality ($P^*=JPJ$), and $1<p<\infty$ by interpolation.
\ef

\begin{theorem}\label{large3} Assume $n=3$ and that (\ref{ass1}) is satisfied with $\rho>1$. One has for $1\le p < q \le \i$
\be \label{e4.10}
\|e^{-tP}\|_{p \to q}\le C t^{-\f{3}{2p}(1-\f p q)}
\ee
for $t\in[1, \i[$.
\end{theorem}

Under the conditions  of Theorem \ref{large3}, it is proved in \cite{w3} that for $s> \f 3 2$,
 \be \label{e4.2w}
 \|\w{x}^{- s} e^{-t P} \w{x}^{- s}\|_{L^2 \to L^2} \le C \w{t}^{- \f 3 2}, t\ge 0.
 \ee
It follows that for $0<r\le \f 3 2$ and $s>r$, one has
 \be \label{e4.3}
 \| e^{-t P} \|_{\vL^{2,s} \to \vL^{2,-s}} \le C \w{t}^{- r}, t\ge 0.
 \ee
 Here $\vL^{2,s} = L^2(\bR^{2n}_{x,v}, \w{x}^{2s}dxdv)$.\\

 {\noindent \bf Proof of Theorem \ref{large3}}
 To obtain large time  $L^p-L^q$ estimate for $e^{-tP}$, we use the following decomposition
 which is deduced from Duhamel's formula :
 \be\label{Duh}
 e^{-tP} = e^{-tP_0} -  I(t) + J(t)
 \ee
 where
 \bea
 I(t) & =& \int_0^t e^{-(t-s)P_0} W e^{-sP_0}\; ds, \label{I}\\
 J(t) &=& \int_0^t \int_0^s e^{-(t-s)P_0} W e^{- \tau P}  W  e^{-(s-\tau) P_0}\; d\tau ds. \label{Jt-def}
 \eea
 Decompose   $I(t) = I_1(t) + I_2(t)$ and $J(t) = J_1(t) + J_2(t)$ where
 \begin{align*}
 I_1(t) &=  \int_0^{\f t 2} e^{-(t-s)P_0} W e^{-sP_0}\; ds,\\
  I_2(t) & =  \int_{\f t 2}^{ t } e^{-(t-s)P_0} W e^{-sP_0}\; ds, \\
 J_1(t) &= \int_0^{\f t 2} \int_0^s e^{-(t-s)P_0} W e^{-\tau P}  W  e^{-(s-\tau) P_0}\; d\tau ds , \\
 J_2(t) &= \int_{\f t 2}^t \int_0^s e^{-(t-s)P_0} W e^{-\tau P}  W  e^{-(s-\tau) P_0}\;  d\tau  ds.
 \end{align*}

We  estimate each term on the right hand side in  $L^1-L^\i$ norm.  Remark first that since $\nabla V(x)$ is bounded,
 the portion of the integral in $I_1(t)$ related to $s\in [0, \f 1 4]$ can be bounded by
 \beas
 \lefteqn{\|\int_0^{\f 1 4}   e^{-(t-s)P_0}  We^{-sP_0} \; ds \|_{1 \to \i}}\\
&\le& \int_0^{\f 1 4} \|\nabla V\|_{\i}   \|e^{-(t-s)P_0}  \|_{1 \to \i } \| \nabla_v e^{-sP_0}\|_{1 \to 1} \; ds\\
  &\le &   Ct^{-\f{3}{2}} \int_0^{\f 1 4} s^{-\f 1 2} \; ds \le  Ct^{-\f{3}{2}}
 \eeas
 for $t \ge 1$. Under the assumption (\ref{ass1}),  $\nabla V(x) \in L^r(\bR^3)$ for any $r > \f{3}{1+\rho}$ and $r\ge 1$.
By H\"older's inequality and (\ref{v3}),  $\nabla V e^{-\delta P_0}$, $\delta>0$,  maps continuously $L^p $ to $L^1$ where
\[
3 <p = \f{r}{r-1} < 1 +  \f{1+\rho}{2-\rho}.
\]
This is possible, because $\f{1+\rho}{2-\rho} > 2$ for $\rho >1$. It follows that
\bea
\| \nabla V {\cdot\na_v} e^{-sP_0}\|_{1 \to 1}   &\le &\| \nabla V e^{-\f 1 8 P_0}\|_{p\to 1}  \| {\na_v} e^{-(s-\f 1 8)P_0}\|_{1 \to p}
\\
&\le &C s^{-\f 3 2 (1 - \f 1 p)}(1+ s^{-\f 1 2}). \nonumber
\eea
This shows that $ s\to \| \nabla V {\cdot\na_v} e^{-sP_0}\|_{1 \to 1} $ is integrable in $s\in [\f 1 4, \i[$.
Consequently, $I_1(t)$ can be bounded as follows
 \beas
\|I_1(t)\|_{1 \to \i}
 & \le &{C}t^{-\f{3}{2}}  + \int_{\f 1 4}^{\f t 2}   \|{ e^{-(t-s)P_0}}  \|_{1 \to \i} \|\nabla V{\cdot\na_v}e^{-sP_0}\|_{1 \to 1} \; ds  \\
& \le &   {C_1}t^{-\f{3}{2}}  \left(1 +   \int_{\f 1 4}^{\f t 2}    \|\nabla V{\cdot\na_v} e^{-sP_0}\|_{1 \to 1} \; ds \right)\\
& \le &  {C_2} t^{-\f{3}{2}}, \quad \mbox{ for } t \ge 1.
 \eeas
Since $\nabla V \in L^r(\bR^{3}_x)$, by (\ref{v4}), $ e^{- \delta P_0}\nabla V$  is bounded from $L^\i$ to $L^r$.
Corollary \ref{cor3.2} shows
\beas
\|e^{- (t-s) P_0}\nabla V\|_{\i \to \i}   &\le  & \|{e^{- (t-s -\delta) P_0}}\|_{r \to \i } \|e^{- \delta P_0}\nabla V \|_{\i \to r} \\
\le C  \w{t-s}^{-\f{3}{2r}} &= & C \w{t-s}^{-\f{1+\rho } 2 + \ep},
\eeas
for  $s \in [\f t 2, t -\f 1 4]$ and $t \ge 1$.
Therefore,  $I_2(t)$ can be estimated by
 \beas
\|I_2(t)\|_{1 \to \i}
 & \le &{C}t^{-\f{3}{2}} + \int^{ t-  \f 1 4}_{\f t 2}   \|e^{-(t-s)P_0} W e^{-s P_0}\|_{1 \to \i} \; ds \\
& \le &   {C_1}t^{-\f{3}{2}}  \left(1 +    \int^{t-\f 1 4}_{\f t 2}   \w{t-s}^{-\f{1+\rho } 2 + \ep}\; ds \right)\\
& \le &  {C_2}t^{-\f{3}{2}}, \quad \mbox{ for } t \ge 1.
 \eeas
It follows that
\be \label{It}
\|I(t)\|_{1 \to \i} \le  Ct^{-\f{3}{2}}, \quad \mbox{ for } t \ge 1.
\ee

For the term $J_1(t)$, we split the domain of integration
$$\Omega = \{( \tau, s );  0 \le s \le \f t 2, 0 \le \tau \le s \} $$ into two parts:  $\Omega = \Omega_1 \cup \Omega_2$, where
\[
\Omega_1 = \{ ( \tau, s) \in \Omega;  \mbox{ either $ \tau \le \f 1 4$ or  $s-\tau \le \f 1 4$}\},
\quad \Omega_2 = \Omega\setminus \Omega_1.
\]
By Corollary \ref{cor3.2} and the fact that $e^{-tP}$ is {uniformly bounded} in $L^p$ {(Proposition \ref{prop:global})}, one can show as above
 that the $L^1-L^\i$ norm of the piece of $J_1(t)$ related to the integration  with respect to  $(\tau,s) \in \Omega_1$ can be bounded by $ C t^{-\f{3}{2}}$.

To treat the remaining part, let $p>3$ be close enough to $3$. Then
\beas
\lefteqn{\|\int_{\Omega_2}  e^{-(t-s)P_0} W e^{-\tau P}  W  e^{-(s-\tau) P_0}\; d\tau ds \|_{1 \to \i} }\\
& \le &  Ct^{-\f{3}{2}} \int_{\Omega_2} \|  e^{- \f 1 8 P_0} \nabla  V e^{- \tau P}  \nabla V e^{-\f 1 8 P_0}\|_{p \to 1}
\w{s-\tau}^{-\f{3(p-1)}{2p}} \; d\tau ds \\
& \le & C_1  t^{-\f{3}{2}} \int_{\f 1 4}^{\f t 2} \| e^{- \f 1 8 P_0}\nabla  V e^{- \tau P}  \nabla V e^{-\f 1 8 P_0}\|_{p \to 1} d\tau
\eeas
 $\w{x}^{-(\f 1 2 + \ep)} : L^p(\bR^3) \to L^2(\bR^3)$ is bounded.
 Condition (\ref{ass1}) and (\ref{v3}) show that
$\nabla V e^{-\f 1 8 P_0}$ is bounded from
$L^p(\bR^6)$ to $\vL^{2, \rho +\f 1 2 -\ep} $. By (\ref{v4}), $e^{-\f 1 8 P_0}\nabla V $ is bounded from $\vL^{2, -\rho + \f 1 2 +\ep}$ to $L^1$. Using (\ref{e4.3}),
we obtain
\[
\| e^{- \f 1 8 P_0}\nabla  V e^{- \tau P}  \nabla V e^{-\f 1 8 P_0}\|_{p \to 1} \le C_\ep\w{\tau}^{-\rho + \f 1 2 +\ep}
\]
It follows that
\be \label{J1}
\|J_1(t)\|_{L^1\to  L^\i}\le \left\{ \begin{array}{ll}
{C}t^{-\f{3}{2}}, & \mbox{ if } \rho > \f 3 2 \\
{C_\ep}t^{ -\rho +\ep}, &   \mbox{ if } 1 < \rho \le \f 3 2.
\end{array}
\right.
\ee
The same estimates hold true for $J_2(t)$. Putting them together, we obtain
\be \label{Jt}
\|J(t)\|_{L^1\to  L^\i}\le \left\{ \begin{array}{ll}
{C}t^{-\f{3}{2}}, & \mbox{ if } \rho > \f 3 2 \\
{C_\ep}t^{ -\rho +\ep}, &   \mbox{ if } 1 < \rho \le \f 3 2.
\end{array}
\right.
\ee
From Corollary \ref{cor3.2}, (\ref{Duh}), (\ref{It}) and (\ref{Jt}), we obtain
\[
\|e^{-tP}\|_{L^1\to  L^\i}\le \left\{ \begin{array}{ll}
{C}t^{-\f{3}{2}}, & \mbox{ if } \rho > \f 3 2 \\
{C_\ep}t^{-\rho +\ep}, &   \mbox{ if } 1 < \rho \le \f 3 2.
\end{array}
\right.
\]
If $\rho > \f 3 2$, then Theorem \ref{large3} is proved by interpolation and the uniform bounds $\|e^{-tP}\|_{p\to p}\le B$ of Proposition \ref{prop:global}.  If $1 <\rho \le \f 3 2$, one obtains
\be \label{beta0}
\|e^{-tP}\|_{L^p\to  L^q}\le
{C_\ep}t^{ - \rho (\f 1 p - \f 1 q)  +\ep},
\ee
for $t >1$ and $1 \le p < q \le \i$.
\\

We now use (\ref{beta0}) instead of (\ref{e4.3}) to improve large-time decay of $e^{-tP}$ for $\rho \in ]1, \f 3 2]$.
Let  $\beta_0(p, q) = \rho(\f 1 p - \f 1 q)$.
Using (\ref{Duh}) and the results for $e^{-tP_0}$, we can show  as before
\be \label{recurrence}
\|e^{-tP}\|_{1\to \i} \le C t^{-\f 3 2} + \|J_1(t)\|_{1 \to \i} + \|J_2(t)\|_{1\to \i}
\ee
 for $t\ge 1$. To estimate $ \|J_1(t)\|_{1 \to \i}$, again we split $\Omega = \Omega_1 \cup \Omega_2$ as before. The $L^1-L^\i$ norm of the piece of $J_1(t)$ given by integration over $\Omega_1$ is bounded by $Ct^{-\f 3 2}$. For $(\tau,s) \in\Omega_2$,  we have
\beas
\lefteqn{\|e^{-(t-s)P_0} W e^{-{\tau} P}  W  e^{-{(s-\tau)} P_0}\|_{1 \to \i}} \\
& \le & \|e^{-(t-s-\delta)P_0}\nabla_v\|_{1\to\i} \|e^{-\delta P_0}\nabla V e^{-{\tau} P}\nabla Ve^{-\delta P_0}\|_{p \to 1} \| \nabla_v e^{-{(s-\tau-\delta)} P_0}\|_{1 \to p},
\eeas
where $p=3+ \ep'$, $\ep'>0$, $\delta>0$. By H\"older's inequality, (\ref{v3}) and (\ref{v4}), $\nabla V e^{-\delta P_0}$ and  $e^{-\delta P_0} \nabla V $ are bounded from $L^p \to L^{1}$, where
\[
\f 1 p + \f 1 {r} =1
\]
and $r >  \f{3}{1+\rho}$.
For $r$ close to $\f{3}{1+\rho}$,
$p$ can be any number smaller than $\f{3}{2-\rho}\in ]3, 6]$ for $1<\rho \le \f 3 2$. Set
\[
r_1 = \beta_0( 1, \f{3}{2-\rho}) =  \f{\rho (1+\rho)}{3}.
\]
 One has $r_1>r_0:=\rho-\f12$.

 Making use of  (\ref{beta0}) (with $(p,q) \sim (1,p)$) instead of (\ref{e4.3}), one obtains
\[
\|J_1(t)\|_{1\to \i} \le \f{C }{t^{\f 3 2}} \left(1+ \int_{\f 1 4}^{\f t 2} s^{-r_1 + \ep} ds\right).
\]
In a similar way,  one can show that $\|J_2(t)\|_{1\to\i}$ satisfies the same estimate.
If  $\rho >\f{-1 +\sqrt{13}}2$,  then $r_1 >1$ and  Theorem \ref{large3} is proved.  If $1 <\rho \le \f{-1 +\sqrt{13}}2$,  one obtains for any $\ep >0$
\be \label{beta1}
\|e^{-tP}\|_{L^p\to  L^q}\le
{C_\ep}t^{ - (\f 1 2 + r_1)(\f 1 p - \f 1 q)  +\ep},
\ee
for $t >1$ and $1 \le p < q \le \i$.
Set $\beta_1(p, q) = (\f 1  2 + r_1)(\f 1 p - \f 1 q) $  and
\[
r_2 = \beta_1( 1, \f{3}{2-\rho}) = (\f 1  2 + r_1)\f{1+\rho}3.
\]
Repeating the arguments from (\ref{beta0}) to (\ref{beta1})  with  (\ref{beta0}) replaced by (\ref{beta1}), one concludes that if $r_2 >1$,
Theorem \ref{large3} is proved. Otherwise, one has
\be \label{beta2}
\|e^{-tP}\|_{L^p\to  L^q}\le
{C_\ep}t^{ - (\f 1 2 +r_2)(\f 1 p - \f 1 q)  +\ep},
\ee
for $t >1$ and $1 \le p < q \le \i$. For $k \ge 3$, set $\beta_{k-1}(p, q) = (\f 1  2 + r_{k-1})(\f 1 p - \f 1 q)$ and
\be \label{rk}
r_k =  \beta_{k-1}(1, \f{3}{2-\rho})= (\f 1  2 + r_{k-1})\f{1+\rho}3 .
\ee
Let $\rho>1$ be fixed. By an induction on $k$, one can prove  that for each $k$, either $r_k >1$, then (\ref{e4.10})
is proved by the above argument; or $0<r_k \le 1$, then one has
\be \label{betak}
\|e^{-tP}\|_{L^p\to  L^q}\le
{C_\ep}t^{ - (\f 1 2 +r_k)(\f 1 p - \f 1 q)  +\ep},
\ee
for $t >1$ and $1 \le p < q \le \i$.
We affirm that for each $\rho>1$, there exists $k\in \bN$ such that $r_k >1$.
In fact, if $r_k \le 1 $ for all $k\in \bN$, then  $\{r_k\}$ would be  an increasing sequence bounded by $1$. Let
$\ell = \lim_{k \to \i} r_k $. Then $\ell \in ]0, 1]$. However, taking the limit $k \to \i$ in (\ref{rk}), one has
\[
\ell = (\f 1  2 + \ell)\f{1+\rho}3,
\]
which gives $ \ell = \f{1+ \rho}{2(2-\rho)}  >1$, because $\rho>1$.
This  contradiction in $\ell$ proves that for each $\rho >1$, there exists some $k$ such that $r_k >1$.
 Therefore (\ref{e4.10}) follows by repeating  at most $k$ times the arguments from (\ref{beta0}) to (\ref{beta1}) with (\ref{beta0}) replaced by newly improved estimate. This achieves the proof of Theorem \ref{large3} for any $\rho>1$.
\ef

Theorem \ref{th1.1} follows from Theorems \ref{short} and \ref{large3}.

\begin{rmk}
Let $n=1$ and condition (\ref{ass1}) be satisfied with $\rho >4$. It is known (\cite{nw}) that for $s> \f 5 2$
\be \label{e3.0}
\|e^{-tP}\|_{{\vL^{2,s} \to \vL^{2,-s}}} \le C\w{t}^{-\f 1 2}, t \ge 1.
\ee
The method used in the proof of Theorem \ref{large3} does not allow to deduce from (\ref{e3.0}) any decay of $e^{-tP}$ in $L^1 -L^\infty$ for $t$ large. For example,
for the  term $I_1(t)$ given in {(\ref{I})},  the method used in the proof of Theorem \ref{large3}  only gives
\[
\|I_1(t)\|_{1 \to\i} \le C \left(t^{-\f 1 2} + \int_{\f 1 4}^{\f t 2} \w{t-s}^{-\f 1 2} \w{s}^{-\f 1 2} ds \right), \quad t \ge 1.
\]
The last integral does not decay as $t\to\i$.
\end{rmk}

\end{document}